\documentclass[twoside,psamsfonts,a4paper]{amsart}
\usepackage{ucs}
\usepackage[utf8x]{inputenc}

\usepackage{amsmath}
\usepackage{amscd}
\usepackage{amsxtra}
\usepackage{amssymb}
\usepackage{amsthm}
\usepackage{stmaryrd}
\usepackage{amsbsy}

\newtheorem{theorem}{Theorem}[section]
\newtheorem{corollary}[theorem]{Corollary}
\newtheorem{lemma}[theorem]{Lemma}
\newtheorem{proposition}[theorem]{Proposition}
\theoremstyle{definition}
\newtheorem{definition}[theorem]{Definition}

\theoremstyle{remark}
\newtheorem{example}[theorem]{Example}
\newtheorem{remark}[theorem]{Remark}

\DeclareSymbolFont{largeletters}{OT1}{cmr}{m}{n}

\DeclareMathOperator{\supp}{supp}
\DeclareMathSymbol{\lcolon}{\mathopen}{operators}{`:}
\DeclareMathSymbol{\rcolon}{\mathclose}{operators}{`:}
\DeclareMathSymbol{\symm}{\mathop}{largeletters}{`S}
\providecommand{\abs}[1]{\lvert#1\rvert}

\newcommand{\Cov}[2]{{#1}{#2}}

\newcommand{\BM}{\mathrm{BM}}
\newcommand{\PD}{\mathrm{PD}}

\newcommand{\ch}{\operatorname{ch}}

\newcommand{\dual}{\spcheck}
\newcommand{\ket}[1]{\lvert#1\rangle}
\newcommand{\bracket}[1]{\langle#1\rangle}
\newcommand{\length}{\ell}
\newcommand{\lie}[1]{\mathfrak{#1}}

\newcommand{\End}{\mathrm{End}}
\newcommand{\Kummer}[2]{{#2}^{[[#1]]}}
\newcommand{\CY}[2]{{#2}^{\{#1\}}}
\newcommand{\Hilb}[2]{{#2}^{[#1]}}
\newcommand{\Set}[1]{\left\{#1\right\}}
\newcommand{\Symm}[2]{{#2}^{(#1)}}
\newcommand{\td}{\operatorname{td}}

\newcommand{\id}{\mathrm{id}}
\newcommand{\pr}{\mathrm{pr}}

\newcommand{\SG}[1]{\mathfrak{S}_{#1}}
\newcommand{\sheaf}[1]{\mathcal #1}
\newcommand{\nop}[1]{\lcolon#1\rcolon}
\newcommand{\scp}[2]{\langle#1, #2\rangle}

\newcommand{\set}[1]{\mathbf{#1}}
\newcommand{\gen}[1]{\boldsymbol{#1}}

\newcommand{\adjoint}{\dagger}

\newcommand{\state}[1]{|#1\rangle}
\newcommand{\vac}{\state{0}}

\title{Twisted cohomology of the Hilbert schemes of points on surfaces}
\author{Marc A.~Nieper-Wi\ss kirchen}
\address{Max-Planck-Insitut for Mathematics, Vivatsgasse 7, 53111 Bonn, Germany}
\email{marc@nieper-wisskirchen.de}
\thanks{The author would like to thank the Max Planck Institute for Mathematics
in Bonn for its hospitality and support during the preparation of this paper.}
\date{\today}

\newcommand{\e}[1]{e_{(#1)}}

\begin{document}
	
	\begin{abstract}
	    We calculate the cohomology spaces of the Hilbert schemes of points on surfaces
	    with values in locally constant systems. For that purpose, we generalise I.~Grojnoswki's
	    and H.~Nakajima's description of the ordinary cohomology in terms of a Fock space
	    representation to the twisted case. We further generalise M.~Lehn's work on
	    the action of the Virasoro algebra to the twisted case.
	    
	    Building on work by M.~Lehn and Ch.~Sorger, we then give an explicit description
	    of the cup-product in the twisted case whenever the surface has a numerically trivial canonical divisor.
	    
	    We formulate our results in a way that they apply to the projective and non-projective
	    case in equal measure.
	    
	    As an application of our methods, we give explicit models for the cohomology rings
	    of the generalised Kummer varieties and of a series of certain even dimensional Calabi--Yau
	    manifolds.
	\end{abstract}
	
	\maketitle
	
	\tableofcontents
	
	\section{Introduction and results}

    Let $X$ be a quasi-projective smooth surface over the complex numbers.
    We denote by $\Hilb n X$ the Hilbert scheme of $n$ points on $X$,
    parametrising zero-dimensional subschemes of $X$ of length
    $n$. It is a quasi-projective variety (\cite{grothendieck:hilbert})
    and smooth of dimension $2n$ (\cite{fogarty}). Recall that the Hilbert scheme
    $\Hilb n X$ can viewed as
    a resolution of the $n$-th symmetric power $\Symm n X := X^n/\SG n$ of the
    surface $X$ by virtue of the Hilbert--Chow morphism $\rho\colon \Hilb n X \to \Symm n X$,
    which maps each zero-dimensional subscheme $\xi$ of $X$ to its support $\supp \xi$ counted with multiplicities.

    Let $L$ be a locally constant system (always over the complex numbers and of rank $1$) over $X$.
    We can view it as a functor from the fundamental groupoid $\Pi$ of $X$ to the
    category of one-dimensional complex vector spaces.
    
    The fundamental groupoid $\Symm n \Pi$ of $\Symm n X$ is the quotient groupoid of $\Pi^n$ by
    the natural $\SG n$-action by \cite{brown} or (in terms of the fundamental group)~\cite{beauville:kummer}.
    Thus we can construct from $L$ a locally constant system $\Symm n L$ on $\Symm n X$ by setting
    \[
        \Symm n L(x_1, \ldots, x_n) := \bigotimes_i L(x_i),
    \]
    for each $(x_1, \ldots, x_n) \in \Symm n X$ (for the notion of the tensor product over an unordered index set
    see, e.g.,~\cite{lehn-sorger:k3}).
    This induces the locally free system $\Hilb n L := \rho^* \Symm n L$ on $\Hilb n X$.
    We are interested in the calculation of the cohomology space
    $\bigoplus_{n \ge 0} H^*(\Hilb n X, \Hilb n L[2n])$. Besided the natural grading given by the cohomological
    degree it carries a weighting given by the number of points $n$. Likewise, the symmetric algebra
    $S^*(\bigoplus_{\nu \ge 1} H^*(X, L^\nu[2]))$ carries a grading by cohomological degree and a weighting, which
    is defined such that $H^*(X, L[2]^\nu)$ is of pure weight $\nu$.
    
    The first result of this paper is the following:
    \begin{theorem}
        \label{thm:vector_space}
    	There is a natural vector space isomorphism
        \[
            \bigoplus_{n \ge 0} H^*(\Hilb n X, \Hilb n L[2n])
            \to
            S^*\left(\bigoplus_{\nu \ge 1} H^*(X, L^\nu[2]\right)
        \]
        that respects the grading and weighting. 
    \end{theorem}
    For $L = \set C$, the trivial system, this result has already
    appeared in~\cite{grojnowski} and~\cite{nakajima}).
    
    Theorem~\ref{thm:vector_space} is proven by defining a Heisenberg Lie algebra $\lie h_{X, L}$,     
    whose underlying vector space
    is given by
    \[
        \bigoplus_{n \ge 0} H^*(X, L^n[2]) \oplus \bigoplus_{n \ge 0} H^*_c(X, L^{-n}[2])
        \oplus \set C \gen c \oplus \set C \gen d
    \]
    and by
    showing that $\bigoplus_{n \ge 0} H^*(\Hilb n X, \Hilb n L[2n])$
    is an irreducible lowest weight representation of this Lie algebra, as is done in~\cite{nakajima}
    for the untwisted case.
    
    Let $G$ be a finite subgroup of the group of locally constant systems on $X$. Via the mapping
    $L \to \Hilb n L$, which is in fact an isomorphism between the groups of locally constant systems on
    $X$ and $\Hilb n X$, respectively,
    $G$ becomes a subgroup of the group of locally constant systems on $\Hilb n X$.
    Such a group naturally defines a
    Galois
    covering $\eta\colon \Cov G {\Hilb n X} \to \Hilb n X$ of degree $\abs G$
    with $\eta_* \set C = \bigoplus_{L \in G} L$. Let us call this covering the $G$-covering of $\Hilb n X$. (In case
    that $G$ is the group of all locally constant systems on $\Hilb n X$, the $G$-covering is the universal one for
    $n \ge 2$.)
    Using the Leray spectral sequence for $\eta$, which already degenerates at the $E_2$-term, the
    cohomology of $\Cov G {\Hilb n X}$ can be computed by Theorem~\ref{thm:vector_space}:
    \begin{corollary}
        \label{cor:vector_space}
    	There is a natural vector space isomorphism
    	\[
            \bigoplus_{n \ge 0} H^*(\Cov G{\Hilb n X}, \set C[2n])
            \to
        	\bigoplus_{L \in G} S^*\left(\bigoplus_{\nu \ge 1} H^*(X, L^\nu[2])\right)
        \]
        that respects the grading and weighting.
    \end{corollary}
    
    We then proceed in the paper
    by defining a twisted version $\lie v_{X, L}$ of the Virasoro Lie algebra, whose underlying
    vector space will be given by
    \[
        \bigoplus_{n \ge 0} H^*(X, L^n) \oplus \bigoplus_{n \ge 0} H^*_c(X, L^{-n})
        \oplus \set C \gen c \oplus \set C \gen d.
    \]
    (Note the different grading compared to $\lie h_{X, L}$.)
    We define an action of $\lie v_{X, L}$
    on $\bigoplus_{n \ge 0} H^*(\Hilb n X, \Hilb n L[2n])$ by generalising results of~\cite{lehn:tautological}
    to the twisted, not necessarily projective case. As in~\cite{lehn:tautological}, we calculate the
    commutators of the operators in $\lie h_{X, L}$ with the boundary operator $\partial$ that
    is given by multiplying
    with $-\frac 1 2$ of the exceptional divisor class of the Hilbert--Chow morphism.
    It turns out that the same relations as in the untwisted, projective case hold.
    
    The next main result of the paper is a decription of the ring structure whenever $X$ has a numerically
    trivial divisor. Following ideas in~\cite{lehn-sorger:k3}, we introduce a family of explicitely
    described graded unital
    algebras $\Hilb n H$ associated to a $G$-weighted (non-counital)
    graded Frobenius algebra $H$ of degree $d$. For example,
    $H = \bigoplus_{L \in G} H^*(X, L[2])$ is such a Frobenius algebra of degree $2$.
    The following holds for each $n \ge 0$:
    \begin{theorem}
        \label{thm:ring}
    	Assume that $X$ has a numerically trivial canonical divisor. Then there is a natural isomorphism
    	\[
    	    \bigoplus_{L \in G} H^*(\Hilb n X, \Hilb n L[2n]) \to \Hilb n {\left(\bigoplus_{L \in G} H^*(X, L[2])\right)}
    	\]
    	of ($G$-weighted) graded algebras of degree $2n$.
    \end{theorem}
    For $L = \set C$, and $X$ projective, this theorem is the main result in~\cite{lehn-sorger:k3}.
    
    The idea of the proof of Theorem~\ref{thm:ring}
    is not to reinvent the wheel but to study how everything can already
    be deduced from the more special case considered in~\cite{lehn-sorger:k3}.
    
    Again by the Leray spectral sequence, also Theorem~\ref{thm:ring}
    has a natural application to the cohomology ring of the
    $G$-coverings of $\Hilb n X$:
    \begin{corollary}
        \label{cor:ring}
    	There is a natural isomorphism
    	\[
    	    H^*(\Cov G {\Hilb n X}, \set C[2n]) \to \Hilb n {\left(\bigoplus_{L \in G} H^*(X, L[2])\right)}
    	\]
    	of graded unital algebras of degree $2n$.
    \end{corollary}

    We want to point out at least two applications of our results. The first one is the computation
    of the cohomology ring of certain families of Calabi--Yau manifolds of even dimension: Let $X$ be an Enriques surface.
    Let $G$ be the
    group of all locally constant systems on $X$, i.e.~$G \simeq \set Z/(2)$. We denote the non-trivial element in $G$ by
    $L$. The Hodge diamonds of $H^*(X, \set C[2])$ and $H^*(X, L[2])$ are given by
    \[
        \begin{matrix}
            & & 1 \\ & 0 & & 0 \\ 0 & & 10 & & 0 \\ & 0 & & 0 \\ & & 1
        \end{matrix}
        \qquad\text{and}\qquad
        \begin{matrix}
            & & 0 \\ & 0 & & 0 \\ 1 & & 10 & & 1, \\ & 0 & & 0 \\ & & 0
        \end{matrix}        
    \]
    respectively.
    
    Denote by $\CY n X$ the (two-fold) universal cover of $\Hilb n X$. By Remark~\ref{rmk:hodge}, the isomorphism
    of Corollary~\ref{cor:vector_space} is in fact an isomorphism of Hodge structures. It follows that
    \[
        H^{k, 0}(\CY n X, \set C) = \begin{cases}
            \set C & \text{for $k = 0$ or $k = 2n$, and}
            \\
            0 & \text{for $0 < k < 2n$.}
        \end{cases}
    \]
    In conjunction with Corollary~\ref{cor:ring}, we have thus proven:
    \begin{proposition}
        \label{prop:cy}   
        For $n > 1$, the manifold $\CY n X$ is a Calabi--Yau manifold in the strict sense. Its cohomology ring
        $H^*(\CY n X, \set C[2n])$ is naturally isomorphic to
        $\Hilb n {(H^*(X, \set C[2]) \oplus H^*(X, L[2]))}$.
    \end{proposition}
    
    Our second application is the calculation of the cohomology ring of the generalised Kummer varieties $\Kummer n X$
    for an abelian surface $X$.
    (A slightly less explicit description of this ring has been obtained by more special methods in~\cite{britze}.)
    Recall from~\cite{beauville:kummer}
    that the generalised Kummer variety $\Kummer n X$ is defined as the fibre over $0$ of the morphism
    $\sigma\colon \Hilb n X \to X$, which is the Hilbert--Chow morphism followed by the summation morphism
    $\Symm n X \to X$ of the abelian surface. The generalised Kummer surface
    is smooth and of dimension $2n - 2$ (\cite{beauville:kummer}).

    As above, let $H$ be a $G$-weighted graded Frobenius algebra of degree $d$. Assume further that $H$ is equipped with
    a compatible structure of a Hopf algebra of degree $d$. For each $n > 0$, we associate to such an algebra an explicitely
    described graded unital algebra $\Kummer n H$ of degree $n$.

    In the following Theorem, we view $H^*(X, \set C[2])$ as such an algebra (the Hopf algebra structure
    given by the group structure of $X$), where we give $H^*(X, \set C[2])$ the trivial
    $G$-weighting for the group $G := X[n]\dual$, the character group of the group of $n$-torsion points on $X$.
    We prove the following:
    \begin{theorem}
        \label{thm:kummer}
        There is a natural isomorphism
        \[
            H^*(\Kummer n X, \set C[2n]) \to \Kummer n {(H^*(X, \set C[2]))}
        \]
        of graded unital algebras of degree $2n$.
    \end{theorem}

    We should remark that most the ``hard work'' that is hidden behind the scenes
    has already been done by others (\cite{grojnowski}, \cite{nakajima},
    \cite{lehn:tautological}, \cite{li-qin-wang:vertex}, \cite{lehn-sorger:k3}, etc.), and our own contribution is to
    see how the ideas and results in the cited papers can be applied an generalised
    to the twisted and to the non-projective case.
    
    \begin{remark}
    	Let us finally mention that the restriction to algebraic, i.e.~quasi-projected surfaces, is unnecessary.
    	Our methods work equally well when we replace $X$ by any complex surface. In this case, the Hilbert schemes
    	become the Douady spaces (\cite{douady}).
    \end{remark}
    
    \section{The Fock space description}
    
    In this section, we prove Theorem~\ref{thm:vector_space} for a locally constant system $L$ on $X$    
    by the method that is used in~\cite{nakajima}
    for the untwisted case, i.e.~by realising the cohomology space of the Hilbert schemes (with coefficients
    in a locally constant system) as an irreducible representation of a Heisenberg Lie algebra.
    
    Let $l \ge 0$ and $n \ge 1$ be two natural numbers. Set
    \[
        \Symm {l, n} X := \Set{(\underline x', x, \underline x) \in \Symm {n + l} X \times X \times \Symm l X\mid
            \underline x' = \underline x + n x}
    \]
    (we write the union of unordered tuples additively).
    The preimage in
    $\Hilb {n + l} X \times X \times \Hilb l X$ of this closed subset under the Hilbert--Chow morphisms $\rho$
    is denoted by
    $\Hilb {n, l} X$ (this incidence variety has already been considered in~\cite{nakajima}).
    
    We denote the projections of $\Symm {l + n} X \times X \times \Symm l X$ onto its three factors by
    $\tilde p$, $\tilde q$ and $\tilde r$, respectively. Likewise, we denote the three projections of
    $\Hilb {l + n} X \times X \times \Hilb l X$ by $p$, $q$ and $r$.
    \begin{lemma}
        \label{lem:trivialisation}
    	It is $q^* L^n \otimes r^* \Hilb l L\vert_{\Hilb {n, l} X} = p^* \Hilb {l + n} L\vert_{\Hilb {n, l} X}$.
    \end{lemma}
    
    \begin{proof}
        It is $\tilde q^* L^n \otimes \tilde r^* \Symm l L\vert_{\Symm {n, l} X}
        = \tilde p^* \Symm {l + n} L\vert_{\Symm {n, l} X}$. This follows from
        \begin{multline*}
            (\tilde q^* L^n \otimes \tilde r^* \Symm l L)(\underline x + n x, x, \underline x)
            \\
            = L(x)^{\otimes n} \otimes \bigotimes_{x' \in \underline x} L(x')
            = \bigotimes_{x' \in \underline x + nx} L(x')
            = \tilde p^* \Symm {l + n} L(\underline x + n x, x, \underline x)
        \end{multline*}
        for every $(\underline x + n x, x, \underline x) \in \Symm {l, n} X$.
        By pulling back everything to the Hilbert schemes, the Lemma follows.
    \end{proof}
    
    Due to Lemma~\ref{lem:trivialisation} and the fact that $p|_{\Hilb {l, n} X}$ is proper (\cite{nakajima}),
    the operator (a correspondence, see~\cite{nakajima})
    \[
        \begin{aligned}
            N\colon H^*(X, L^n[2]) \times H^*(\Hilb l X, \Hilb l L[2l])
            & \to H^*(\Hilb {l + n} X, \Hilb {l + n} L[2 (l + n)]),\
            \\
            (\alpha, \beta) & \mapsto \PD^{-1} p_*((q^* \alpha \cup r^* \beta) \cap [\Hilb {l, n} X])
        \end{aligned}
    \]
    is well defined. Here,
    \[
        \PD\colon H^*(\Hilb {l + n} X, \Hilb {n + l} L[2(l + n)]) \to
        H^\BM_*(\Hilb {l + n} X, \Hilb {n + l} L[-2 (l + n)])
    \]
    is the Poincaré-duality isomorphism between the cohomology
    and the Borel--Moore homology.
    (The degree shifts are chosen in a way that $N$ an operator of degree $0$, see~\cite{lehn-sorger:k3}.)
    
    Furthermore, $q \times r|_{X \times \Hilb l X}$ is proper (\cite{nakajima}). Thus we can also define
    an operator the other way round:
    \[
        \begin{aligned}
            N^\adjoint\colon H^*_c(X, L^{-n}[2]) \times H^*(\Hilb {n + l} X, \Hilb {l + n} L [2])
            & \to H^*(\Hilb l X, \Hilb l L[2l]),
            \\
            (\alpha, \beta) & \mapsto (-1)^n \PD^{-1} r_* (q^* \alpha \cup p^* \beta \cap [\Hilb {l, n} X])
        \end{aligned}
    \]
    
    As in~\cite{nakajima}, we will use these operators to define an action of a Heisenberg Lie algebra on
    \[
        V_{X, L} := \bigoplus_{n \ge 0} H^*(\Hilb n X, \Hilb n L).
    \]
    For this, let $A$ be a weighted, graded
    Frobenius algebra of degree $d$
    (over the complex numbers), that
    is a weighted and graded vector space over $\set C$
    with a (graded) commutative and associative multiplication of degree $d$
    and weight $0$ and a unit element $1$ (necessarily of degree $-d$ and weight $0$) together with a linear form
    $\int\colon A \to \set C$ of degree $-d$ and weight $0$ such that for each weight $\nu \in \set Z$ the induced
    bilinear form $\scp{\cdot}{\cdot}\colon A(\nu) \times A(-\nu) \to \set C, (a, a') \mapsto \int_A a a'$ is
    non-degenerate (of degree $0$). Here $A(\nu)$ denotes the weight space of weight $\nu$. In particular, all weight spaces
    are finite-dimensional. In the case of a trivial weighting, this notion of a graded Frobenius algebra
    has already appeared
    in~\cite{lehn-sorger:k3}.
    \begin{example}
    	The vector space
    	\[
    	    A_{X, L} := \bigoplus_{\nu \ge 0} H^*(X, L^\nu[2]) \oplus \bigoplus_{\nu \ge 0} H^*_c(X, L^{-\nu}[2])
    	\]
    	is naturally a weighted, graded Frobenius algebra of degree $2$ ($= \dim X$)
    	as follows: The grading is given by the cohomology grading.
    	The weighting is defined by defining $H^*(X, L^\nu[2])$ of pure weight $\nu$ for $\nu \ge 0$ and
    	$H^*_c(X, L^{-\nu}[2])$ of pure weight $-\nu$. The multiplication is given by the cup-product,
    	where we view the product of an ordinary cohomology class and of a cohomology class with compact support
    	as an ordinary cohomology class whenever the resulting weight is strictly positive,
    	and as a cohomology class with compact support otherwise. The linear form $\int$ is given by
    	evaluating a class with compact support on the fundamental class of $X$.
    \end{example}
    For such a weighted, graded Frobenius algebra $A$ we set
    \[
        \lie h_A := A \oplus \set C \gen c \oplus \set C \gen d.
    \]
    We define the structure of a weighted,
    graded Lie algebra on $\lie h_A$ by defining $\gen c$ to be a central element of weight $0$
    and degree $0$, $\gen d$ an element of weight $0$ and degree $0$ and by setting the following commutator relations:
    $[\gen d, a] := n \cdot a$ for each element $a \in A$ of weight $n$,
    and $[a, a'] = \scp{[\gen d, a]}{a'} \gen c$ for elements $a, a' \in A$.
    \begin{definition}
        The Lie algebra $\lie h_A$ the \emph{Heisenberg algebra associated to $A$}.     	
    \end{definition}
    For $A = A_{X, L}$, we set $\lie h_{X, L} := \lie h_A$. We define a linear map
    \[
        q\colon \lie h_{X, L} \to \End(V_{X, L})
    \]
    as follows:
    Let $l \ge 0$ and $\beta \in V_{X, L}(l) = H^*(\Hilb l X, \Hilb l L[2l])$. We set $q(\gen c)(\beta) := \beta$, and 
    $q(\gen d)(\beta) := l \beta$. For $n \ge 0$,
    and $\alpha \in A_{X, L}(\nu) = H^*(X, L^{\nu}[2])$, we set $q(\alpha)(\beta) := N(\alpha, \beta)$. For
    $\alpha \in A_{X, L}(-\nu) = H^*_c(X, L^{-\nu}[2])$, we set $q(\alpha)(\beta) := N^\adjoint(\alpha, \beta)$.
    Finally, we set $q(\alpha)(\beta) = 0$ for $\alpha \in A_{X, L}(0) = H^*(X, \set C) \oplus H^*_c(X, \set C)$.
    
    \begin{proposition}
        \label{prop:vector_space}
    	The map $q$ is a weighted, graded action of $\lie h_{X, L}$ on $V_{X, L}$.
    \end{proposition}
    
    \begin{proof}
        This Proposition is proven in~\cite{nakajima} for the untwisted case, i.e.~for $L = \set C$. The proof
        there is based on calculating commutators on the level of cycles
        of the correspondences defined by the incidence schemes
        $\Hilb {l, n} X$. These commutators are independent of the locally constant system used. Thus the proof
        in~\cite{nakajima} also applies to this more general case.
    \end{proof}
    
    \begin{example}
        \label{ex:hilbert-chow}
        Let $\alpha = \sum \alpha_{(1)} \otimes \cdots \otimes \alpha_{(n)} \in H^*(\Symm n X, \Symm n L[2n]) = 
        S^n H^*(X, L[2])$ (we use the Sweedler notation to denote elements in tensor products). The pull-back
        of $\alpha$ by the Hilbert--Chow morphism $\rho\colon \Hilb n X \to \Symm n X$ is then given by
        \[
            \rho^* \alpha = \frac 1 {n!} \sum q(\alpha_{(1)}) \cdots q(\alpha_{(n)}) \vac,
        \]
        where $\vac$ is the unit $1 \in H^*(\Hilb 0 X, \set C) = \set C$.
    \end{example}

    We will use Proposition~\ref{prop:vector_space}
    to prove our first Theorem.
    \begin{proof}[Proof of Theorem~\ref{thm:vector_space}]
        The vector space $\tilde V_{X, L} := S^*(\bigoplus_{\nu \ge 1} H^*(X, L^\nu[2]))$
        carries a unique structure as an $\lie h_{X, L}$-module such that $\gen c$ acts as the identity, $\gen d$ acts
        by multiplying with the weight, $\alpha \in H^*(X, L^n)$ for $n \ge 1$ acts by multiplying with $\alpha$,
        and $\alpha \in H^*(X, \set C) \oplus H^*_c(X, \set C)$ acts by zero. By the representation theory of the
        Lie algebras of Heisenberg-type, this is an irreducible lowest weight representation of $\lie h_{V, L}$,
        which is generated by the lowest weight vector $1$, which is of weight $0$.
        
        The $\lie h_{V, L}$-module $V_{X, L}$ also has a vector of weight $0$,
        namely $\vac$.
        Thus, there is a unique morphism $\Phi\colon \tilde V_L \to V_L$
        of $\set h_L$-modules that maps $1$ to $\vac$. This will be the inverse of the 
        isomorphism mentioned in Theorem~\ref{thm:vector_space}.
        It remains to show that $\Phi$ is bijective.
        The injectivity follows from the fact that $\tilde V_{X, L}$ is irreducible as an $\lie h_{X, L}$-module.
        
        In order to prove the surjectivity, we will derive upper bounds on the dimensions of the weight spaces of the
        right hand side $V_{X, L}$ (see also~\cite{lehn:hilbert_lectures} about this proof method).
        By the Leray spectral sequence associated
        to the Hilbert--Chow morphism $\rho\colon \Hilb n X \to \Symm n X$, such an upper bound is provided by
        the dimension of the spectral sequence's $E_2$-term
    	$H^*(\Symm n X, \set R^* \rho_* L[2 n])$. As shown in~\cite{goettsche-soergel}, it follows from the
	    Beilinson--Bernstein--Deligne--Gabber decomposition theorem that
    	\[
    	    \set R^* \rho_* \set Q[2 n] = \bigoplus_{\lambda \in \set P(n)} (i_\lambda)_* \set Q[2\length(\lambda)].
    	\]
    	Here, $P(n)$ is the set of all partitions of $n$, $\length(\lambda) = r$ is the length of a partition
    	$\lambda = (\lambda^1, \lambda^2, \ldots, \lambda^r)$, $\Symm \lambda X := \Set{\sum_{i = 1}^r \lambda_i x_i \mid
    	x_i \in X} \subset \Symm n X$, and $i_\lambda\colon \Symm \lambda X \to \Symm n X$ is the inclusion map.
    	
    	Set $\Symm \lambda L := i_\lambda^* \Symm n L$. By the projection formula, it follows that
    	$\set R^* \rho_*  L[2 n] = \bigoplus_{\lambda \in \set P(n)}
    	(i_\lambda)_* \Symm \lambda L[2 \length(\lambda)]$.
    	
    	Thus, an upper bound on the dimension of $H^*(\Hilb n X, \Hilb n L[2n])$ is provided by the dimension of
    	$\bigoplus_{\lambda \in \set P(n)} H^*(\Symm \lambda X, \Symm \lambda L[2 \length(\lambda)])$.
    	By~\cite{goettsche-soergel}, this can be seen to be isomorphic to
    	\[
    	    \bigoplus_{\sum_{i \ge 1} i \nu_i = n} \bigotimes_{i \ge 1} S^{\nu_i} H^*(X, L^i[2]),
    	\]
    	where each $\nu_i \ge 0$.
    	It follows that the upper bound given by the $E_2$-term is exactly the dimension of
    	the $n$-th weight space of $\tilde V_{X, L}$.
    	Thus the dimension of the weight spaces
    	of $V_{X, L}$ cannot be greater than the dimensions of the weight spaces of $\tilde V_{X, L}$.
    	Thus the Theorem is proven.
    \end{proof}
    
    \begin{remark}
        \label{rmk:hodge}
        Assume that $X$ is projective. In this case, the (twisted)
        cohomology spaces of $X$ and its Hilbert schemes $\Hilb n X$ carry pure Hodge structures. As the isomorphism
        of Theorem~\ref{thm:vector_space} is defined by algebraic correspondences (i.e.~by correspondences of Hodge type
        $(p, p)$), it follows that the isomorphism in Theorem~\ref{thm:vector_space} is compatible with the
        natural Hodge structures on both sides.
        
        In terms of Hodge numbers, the following equation encodes our result:
        \begin{multline*}
            \sum_{n \ge 0} \prod_{i, j} h^{i, j}(\Hilb n X, \Hilb n L[2n]) p^i q^j z^n
            \\
            = \prod_{m \ge 1} \prod_{i, j} (1 - (-1)^{i + j} p^i q^j z^m)^{- (-1)^{i + j} h^{i, j}(X, L^m[2])}
        \end{multline*}
    \end{remark}
    
    \section{The Virasoro algebra in the twisted case}
    
    To each weighted, graded Frobenius algebra $A$ of degree $d$, we associate a skew-symmetric 
    form
    $e\colon A \times A \to \set C$ of degree $d$ as follows:
    
    Let $n \in \set Z$. We note that
    $A(n)$ and $A(-n)$ are dual to each other via the linear form $\int$. Thus we can consider the linear map
    $\Delta(n)\colon \set C \to A(n) \otimes A(-n)$
    dual to the bilinear form $\scp \cdot\cdot\colon A(n) \otimes A(-n) \to \set C$.
    Write $\Delta(n) 1 = \sum e_{(1)}(n) \otimes e_{(2)}(n)$ in Sweedler notation. Then we define $e$
    by setting
    \[
        e(\alpha, \beta) := \sum_{\nu = 0}^n \frac{\nu (n - \nu)} 2 \int \sum e_{(1)}(\nu) e_{(2)}(\nu) \alpha \beta
    \]
    for all $\alpha \in A(n)$ whenever $n \ge 0$. We shall call this form the \emph{Euler form of $A$}.
    \begin{example}
    	Assume that $A(n) \equiv A(0)$ for all $n \in \set Z$. In this case, we have
    	\[
    	    e(\alpha, \beta) = \frac{n^3 - n}{12} \int e\alpha \beta
    	\]
    	for $\alpha \in A(n)$ with
   	    $e := \int \sum e_{(1)}(0) e_{(2)}(0)$ (\cite{lehn:tautological}).
    \end{example}
    
    We use the Euler form to define another Lie algebra associated to $A$. We set
    \[
        \lie v_A := A[-2] \oplus \set C \gen c \oplus \set C \gen d.
    \]
    We define the structure of a weighted, graded
    Lie algebra on $\lie v_A$ be defining $\gen c$ to be a central element or weight $0$ and degree $0$,
    $\gen d$ an alement of
    weight $0$ and degree $0$ and by introducing the following commutator relations:
    $[\gen d, a] := n \cdot a$ for each element $a \in A[-2]$ of weight $n$, and
    $[a, a'] := (\gen d a) a' - a (\gen d a') - e(a, a')$ for elements $a, a' \in A$.
    \begin{definition}
    	The Lie algebra $\lie v_A$ is the \emph{Virasoro algebra associated to $A$}.
    \end{definition}
    For $A = A_{X, L}$, we set $\lie v_{X, L} := \lie v_A$.
    The whole construction is a generalisation to the twisted case
    of the Virasoro algebra found in~\cite{lehn:tautological}.

    We now define a linear map $L\colon \lie v_{X, L} \to \End(V_{X, L})$ as follows:
    We define $L(\gen c)$ to be the identity, $L(\gen d)$ to be multiplication with the weight, and
    for $\alpha \in A[-2]$ we set
    \[L(\alpha) := \frac 1 2 \sum \sum_{\nu \in \set Z}
    \nop{q(e_{(1)}(\nu)) q(e_{(2)}(\nu) \alpha)},\]
    where the \emph{normal
    ordered product $\nop{a a'}$} of two operators is defined to be $a a'$ if the weight of $a$ is greater or equal
    to the weight of $a'$ and is defined to be $a' a$ if the weight of $a'$ is greater than the weight of $a$.
    
    The following Lemma is proven for the untwisted case in~\cite{lehn:tautological}.
    \begin{lemma}
        \label{lem:vir_and_q}
        For $\alpha \in A_{V, L}[-2]$ and $\beta \in A_{V, L}$, we have
        \[
            [L(\alpha), q(\beta)] = - q(\alpha [\gen d, \beta]).
        \]
    \end{lemma}
    
    \begin{proof}
        Let $\alpha \in A_{V, L}[2](n)$ and $\beta \in A_{V, L}(m)$ with $n, m \in \set Z$. In the following
        calculations we omit all Koszul signs arising from commuting the graded elements $\alpha$ and $\beta$.
        By definition,
        we have
        $[L(\alpha), q(\beta)] = \frac 1 2 \sum \sum_{\nu} [\nop{q(e_{(1)}(\nu)) q(e_{(2)}(\nu) \alpha)}, q(\beta)]$,
        where $\nu$ runs through all integers.
        As the commutator of two operators in $\lie h_{V, L}$ is central, we do not have to pay attention
        to the order of the factors of the normally ordered product when calculating the commutator:
        \begin{multline*}
            [\nop{q(e_{(1)}(\nu)) q(e_{(2)}(\nu) \alpha)}, q(\beta)]
            \\
            =  \nu \scp{e_{(1)}(\nu)}{\beta} q(e_{(2)}(\nu) \alpha)
            	+ (n - \nu) \scp{e_{(2)}(\nu) \alpha}{\beta} q(e_{(1)}(\nu)).
        \end{multline*}
        As $\scp\cdot\cdot$ is of weight zero, the first summand is only non-zero for $\nu = -m$, while
        the second summand is only non-zero for $\nu = n + m$. Thus we have
        \begin{multline*}
        	[L(\alpha), q(\beta)] = - \frac m 2 \sum \left(\scp{e_{(1)}(-m)}{\beta} q(e_{(2)}(-m ) \alpha)\right.
            	\\
            	\left. + \scp{e_{(2)}(n + m) \alpha}{\beta} q(e_{(1)}(n + m))\right).
        \end{multline*}
        As $e_{(1)}(\cdot)$ is the dual basis to $e_{(2)}(\cdot)$,
        the right hand side simplifies to $-m q(\alpha \beta)$, which proves the Lemma.
    \end{proof}
    
    We use Lemma~\ref{lem:vir_and_q}
    to prove the following Proposition, which has already appeared in~\cite{lehn:tautological} for
    the untwisted, projective case:
    \begin{proposition}
        The map $L$ is a weighted, graded action of the Virasoro algebra $\lie v_{X, L}$ on $V_{X, L}$.        	
    \end{proposition}
    
    \begin{proof}
    	Let $\alpha \in A[-2](m)$ and $\beta \in A[-2](n)$ with $m, n \in \set Z$. We have to prove that
        $[L(\alpha), L(\beta)] = (m - n) L(\alpha \beta) - e(\alpha, \beta)$.
        We follow ideas in~\cite{frenkel:voas}. In all summations below, $\nu$ runs through all integers if not
        specified otherwise.
    	
    	We begin with the case $n \neq 0$ and $m + n \neq 0$. In this case, by Lemma~\ref{lem:vir_and_q}, it is
    	\begin{multline*}	
    	    	[L(\alpha), L(\beta)]
    	    	= \frac 1 2 \left[L(m), \sum \sum_\nu q(\e1(\nu)) q(\e2(\nu) \beta)\right]
    	    	\\
    	    	= \frac 1 2 \left(\sum \sum_\nu 
    	    	    (-\nu) q(\e1(\nu) \alpha) q(\e2(\nu)\beta)
    	    	    + \sum \sum_\nu (\nu - n) q(\e1(\nu)) q(\e2(\nu) \alpha \beta)\right).
    	\end{multline*}    	
        As $\sum q(\e1(\nu)(\alpha) q(\e2(\nu)\beta) = q(\e1(\nu + m)) q(\e2(\nu - m) \alpha \beta)$,
        the right hand side is equal to
        \[
            \frac 1 2 \sum \sum_\nu \left(
    	    	    (-\nu) q(\e1(\nu + m)) q(\e2(\nu + m)\alpha\beta)
    	    	    + (\nu - n) q(\e1(\nu)) q(\e2(\nu) \alpha \beta)
    	    	\right),
    	\]
    	which is nothing else than $(m - n) L(\alpha \beta)$. Note that $e(\alpha, \beta) = 0$ in this case.
    	
        The next case we study is $m > 0$ and $n = - m$. In order to ensure convergence in the following calculations
        we have to split up $L(\beta)$
        as follows:
        \[
            L(\beta) = \sum \sum_{\nu \ge m} q(\e1(\nu) \beta) q(\e2(\nu)) +
            	    \sum \sum_{\nu < m} q(\e2(\nu)) q(\e1(\nu) \beta)
        \]
        Calculating the commutator $[L(\alpha), L(\beta)]$ thus yields the four terms:
        \begin{multline*}        	
             \frac 1 2 \sum \sum_{\nu \ge m} (m - \nu) q(\e1(\nu) \alpha \beta) q(\e2(\nu))
            	+ \frac 1 2 \sum \sum_{\nu \ge m} \nu q(\e1(\nu) \beta) q(\e1(\nu) \alpha)
            	\\
            	+  \frac 1 2 \sum \sum_{\nu < m} \nu q(\e2(\nu) \alpha) q(\e1(\nu) \beta)
                + \frac 1 2 \sum \sum_{\nu < m} (m - \nu) q(\e2(\nu)) q(\e1(\nu) \alpha\beta).
        \end{multline*}
        As in the first case, we now move $\alpha$ and $\beta$ rightwards. Then we can split off
        an infinite part given by a multiple of $L(\alpha\beta)$ and are left over with the finite sum
        \begin{multline*}
        	[L(\alpha), L(\beta)] - 2m L(\alpha\beta)
        	\\
        	= \frac 1 2 \sum \sum_{\nu = 0}^m (m - \nu) \left(q(\e2(\nu)) q(\e1(\nu) \alpha \beta)
            - q(\e1(\nu)) q(\e2(\nu) \alpha \beta)\right).
        \end{multline*}
        The right side is exactly $e(\alpha, \beta)$.
                
        The remaining cases either follow from the above by exchanging $n$ and $m$ or are trivial ($n = m = 0$).
    \end{proof}

    \section{The boundary operator}
    
    We proceed as in~\cite{lehn:tautological} by introducing a boundary operator on $V_{X, L}$. 
    Recall the definition of the tautological classes of the Hilbert scheme~\cite{li-qin-wang:vertex}:
    Let $\Xi^n$ be the universal family over $\Hilb n X$, which is a subscheme of $\Hilb n X \times X$.
    We denote the projections of $\Hilb n X \times X$ onto its factors by $p$ and $q$.
    To each $\alpha \in H^*(X, \set C)$ we associate the \emph{tautological classes}
    \[
        \Hilb n \alpha := p_*(\ch(\sheaf O_{\Xi^n}) \cup q^*(\td(X) \cup \alpha))
    \]
    in $H^*(\Hilb n X, \set C)$.
    \begin{remark}
    	Note that the tautological classes live in the cohomology with untwisted coefficients, and we have not
    	generalised this concept to the twisted case.
    \end{remark}
    Each $\alpha \in H^*(X, \set C)$ defines an operator $m(\alpha) \in \End(V_{X, L})$,
    which is given by $m(\alpha)(\beta) := \Hilb n \alpha \cup \beta$ for all $\beta \in H^*(\Hilb n X, \Hilb n L)$.
    It is an operator of weight zero. As it does not respect the grading, we split it up into its homogeneous
    components $m(\alpha) = \sum m^*(\alpha)$ with respect to the grading. Following~\cite{lehn:tautological}, we
    set $\partial := m^2(1)$ and call it the \emph{boundary operator}. It is an operator of weight $0$ and degree $2$.
    For each operator $p \in \End(V_{X, L})$,
    we set $p' := [\partial, p]$ and call it the \emph{derivative of $p$}.

    The main theorem in~\cite{lehn:tautological} is the calculation of the derivatives of the Heisenberg operators
    in the untwisted, projective case.
    In the sequel, we will do this in our more general case:

    Let $K$ be the canonical divisor of $X$. We make it into an operator
    $K\colon A_{X, L} \to A_{X, L}[-2]$ of weight zero
    by setting
    \[
        K(\alpha) := \frac{\abs n - 1} 2 K \alpha
    \] for $\alpha \in H^*(\Hilb n X, \Hilb n L)$.
    \begin{proposition}
        \label{prop:lehn_main}
    	For all $\alpha, \beta \in A_{X, L}$ the following holds:
    	\[
    	    [q'(\alpha), q(\beta)] = - q([\gen d, \alpha] [\gen d, \beta])
    	        - \int K([\gen d, \alpha]) [\gen d, \beta].
    	\]
    \end{proposition}
    
    \begin{proof}
    	Let us first consider the case of $\alpha \in A(m)$ and $\beta \in A(n)$ with $n + m \neq 0$. We have to
    	show that $[q'(\alpha), q(\beta)] = -n m q(\alpha \beta)$. This is proven in~\cite{lehn:tautological} for the      
    	projective, untwisted case. The proof in~\cite{lehn:tautological} is based on
    	calculating the commutator on the level of cycles. As these calculations are local in $X$, the result
    	remains true for non-projective $X$. Furthermore, the proof literally works in the twisted case.
    	
    	The case $n + m = 0$ remains. Here we have to show that $[q'(\alpha), q(\beta)] = 
    	m^2 \frac {\abs m - 1}{2} \int K \alpha \beta$. In~\cite{lehn:tautological} the following intermediate
    	result is formulated for the projective, untwisted case: For all $m \in \set Z$, there exists a class
    	$K_m \in H^*(X, \set C)$ such that $[q'(\alpha), q(\beta)] = m^2 \id \int K_m \alpha \beta$.
    	As above the proof for this intermediate result that is given in~\cite{lehn:tautological}
    	also works in the twisted and non-projective case. The classes $K_m$ do not depend
    	on the choice of $L$, i.e.~are universal for the surface. In~\cite{lehn:tautological}, the classes $K_m$ are 
    	computed for the projective case, namely $K_m = \frac{\abs m - 1}2 K$, where $K$ is the class of the canonical
    	divisor. All that remains is to calculate the classes $K_m$ for the non-projective (untwisted) case.
    	As $[q'(\alpha), q(\beta)] = [q'(\beta), q(\alpha)]$ (up to Koszul signs), it is enough to calculate $K_m$ for
    	$m > 0$:
    	
    	Let
    	$\beta \in A_{X, \set C}(-m) = H_c^*(X, \set C[2])$. Consider an open embedding $j\colon X \to \hat X$
    	of $X$ into a smooth, projective surface $\hat X$. We denote
    	the corresponding embeddings $\Hilb n X \to \Hilb n {\hat X}$
    	also by the letter $j$. Denote the $1$ in $A_{\hat X, \set C}(m) = H^*(X, \set C[2])$ by $1(m)$.
    	As all constructions considered so far are functorial (in
    	the appropriate senses) with respect to open embeddings, we have
        \[
            j^* [q'(1(m)), q(j_* \beta)] \vac = [q'(j^* 1(m)), q(\beta)] \vac.
        \]
        The right hand side is given by
        $m^2 \int K_m \beta$, where $K_m$ is the class corresponding to $X$.
        By the calculations in~\cite{lehn:tautological},
        the left hand side is given by $m^2 \frac{\abs m - 1} 2 \int K_{\hat X} j_* \beta$,
        where $K_{\hat X}$ is the canonical divisor class of $\hat X$.
        As $j^*K_{\hat X} = K_X$, we see that $K_m = \frac{\abs m - 1}2 K$ also holds in the non-projective case,
        which proves the Proposition.
    \end{proof}
    
    \begin{corollary}
    	For all $\alpha \in A_{X, L}$, the following holds:
    	\[
    	    q'(\alpha) = L([\gen, d \alpha]) + q(K([\gen, d \alpha])).
    	\]
    \end{corollary}
    
    \begin{proof}
    	This can be deduced from~\ref{prop:lehn_main} as the respective statement for the untwisted, projective case
    	in~\cite{lehn:tautological} is proven.
    \end{proof}
    
    \section{The ring structure}
    
    From now on, we assume that the canonical divisor of $X$ is numerically trivial.
    
    Let $H$ be a \emph{non-counital graded Frobenius algebra of degree $d$ (over the complex numbers)}, that is a
    graded vector space over $\set C$ with a (graded) commutative and associative multiplication of degree
    $d$ and a unit element $1$ (of degree $-d$) together with a coassociative and cocommutative
    $H$-module homomorphism $\Delta\colon H \to H \otimes H$ of degree $d$. (We regard $H \otimes H$ as an $H$-algebra
    by multiplying on the left factor.)
    The map $\Delta$ is called the \emph{diagonal}.
    
    \begin{example}
    	Let $H$ be a graded Frobenius algebra of degree $d$. The dual
    	$\Delta$ to the multiplication map $H \otimes H \to H$
    	with respect to $\scp\cdot\cdot$ makes $H$ a non-counital, graded Frobenius algebra of degree $d$.
    	In this context, the integral $\int$ of the Frobenius algebra is the \emph{counit of $H$}.
    \end{example}

    Let $G$ be a finite abelian group, which will be written additively in the sequel.
    A $G$-weighting on $H$ is an action of the character group $G\dual$ of $G$ on $H$.
    In other words, $H$ comes together with a weight decomposition of the form $\bigoplus_{L \in G} H(L)$,
    where each $\chi \in G\dual$ acts on $H(L)$ by multiplication with $\chi(L)$. 
    \begin{example}
        Let $G$ be a finite subgroup of the group of locally constant systems on $X$, written additively.
    	The $G$-weighted vector space
    	\[
    	    H_{X, G} := \bigoplus_{L \in G} H^*(X, L[2])
    	\]
    	is naturally a non-counital
    	$G$-weighted, graded Frobenius  algebra of degree
    	$2$ as follows: the grading is given by the cohomological grading. The multiplication is given by the cup
    	product. The diagonal is given by the proper push-forward $\delta_*\colon H_{X, G} \to H_{X, G} \otimes H_{X, G}$
    	that is induced by the diagonal map $\delta\colon X \times X \to X$.
    \end{example}
    
    By iterated application, 
    $\Delta$ induce maps $\Delta\colon H \to H^{\otimes n}$ with $n \ge 1$. We denote the restriction
    of $\Delta\colon H \to H^{\otimes n}$ to $H(n L)$, $L \in G$, followed by the projection onto $H(L)^{\otimes n}$
    by $\Delta(L)\colon H(n L) \to H(L)^{\otimes n}$.
    The element $e := (\nabla \circ \Delta(0))(1) \in H$ is called the \emph{Euler class of $H$}, where
    $\nabla\colon H \otimes H \to H$ is the multiplication map.

    There is a construction given in~\cite{lehn-sorger:k3} that associates to each graded Frobenius algebra $H$ of degree
    of $d$ a sequence of graded Frobenius algebras $\Hilb n H$ (whose degrees are given by $nd$). We extend this
    construction to $G$-weighted not necessarily counital Frobenius algebras as follows:
    For each $L \in G$, set
    \[
        H_n(L) := \bigoplus_{\sigma \in \SG n} \left(\bigotimes_{B \in \sigma\backslash[n]} H(\abs B L)\right)
        \gen \sigma\quad\text{and}\quad
        H_n := \bigoplus_{L \in G} H_n(L),
    \]
    where $[n] := \Set{1, \ldots, n}$ and $\sigma\backslash[n]$ is the set of orbits of the action of the cyclic group
    generated by $\sigma$ on the set $[n]$.
    (Note that $H_n(0) = H(0)\{\SG n\}$ in the terminology of~\cite{lehn-sorger:k3}.)
    The symmetric group $\SG n$ acts on $H_n$.
    The graded vector space of invariants, $H_n^{\SG n}$, is denoted by $\Hilb n H$.
    
    Let $f\colon I \to J$ a surjection of finite sets
    and $(n_i)_{i \in I}$ a tuple of integers.
    Fibre-wise multiplication yields ring homomorphisms
    \[
        \nabla^{I, J} := \nabla^f\colon \bigotimes_{i \in I} H(n_i L)
        \to \bigotimes_{j \in J} H\left(\sum_{f(i) = j} n_i L \right)
    \]
    of degree $d(\abs I - \abs J)$. (These correspond to the ring homomorphism $f^{I, J}$ in~\cite{lehn-sorger:k3}.)
    Dually, by using the diagonal morphisms $\Delta(L)$ and relying on their coassociativity and cocommutativity,
    we can define $\nabla^f$-module homomorphisms
    \[
        \Delta_{J, I} := \Delta_f\colon \bigotimes_{j \in J} H\left(\sum_{f(i) = j} n_i L\right)
        \to \bigotimes_{i\in I} H(n_i L),
    \]
    which are also of degree $d(\abs I - \abs J)$. (These correspond to the
    module homomorphisms $f_{J, I}$ in~\cite{lehn-sorger:k3}).
    
    Let $\sigma, \tau \in \SG n$ be two permutations. By $\langle\sigma, \tau\rangle$ we denote the subgroup of $\SG n$
    generated by the two permutations. Note that there are natural surjections
    $\sigma\backslash[n] \to \langle\sigma, \tau\rangle\backslash[n]$,
    $\tau\backslash[n] \to \langle\sigma, \tau\rangle\backslash[n]$, and
    $(\sigma\tau)\backslash[n] \to \langle\sigma, \tau\rangle\backslash[n]$. The corresponding ring and module homomorphism
    are denoted by $\nabla^{\sigma, \langle\sigma, \tau\rangle}$, etc., and $\Delta_{\langle\sigma, \tau\rangle, \sigma}$,
    etc.
    
    Let $L, M \in G$. We define a linear map
    \[
        m_{\sigma, \tau}\colon \bigotimes_{B \in \sigma\backslash[n]} H(\abs B L) \otimes 
        \bigotimes_{B \in \tau\backslash[n]} H(\abs B M) \to
        \bigotimes_{B \in (\sigma\tau)\backslash[n]} H(\abs B (L + M))
    \]
    by
    \[
        m_{\sigma, \tau}(\alpha \otimes \beta) = \Delta_{\langle\sigma, \tau\rangle, (\sigma\tau)}(\nabla^{\sigma,
        \langle\sigma, \tau\rangle}(\alpha) \nabla^{\tau, \langle\sigma, \tau\rangle}(\beta)
        e^{\gamma(\sigma, \tau)}),
    \]
    where the expression $e^{\gamma(\sigma, \tau)}$
    is defined as in~\cite{lehn-sorger:k3} (we have to use our Euler class $e$, which is defined above).
    This defines a product $H_n \otimes H_n \to H_n$ which is given by
    \[
        (\alpha \gen \sigma) \cdot (\beta \gen \sigma) := m_{\sigma, \tau}(\alpha, \beta) \gen {\sigma\tau}
    \]
    for $\alpha \gen \sigma \in H_n(L)$ and $\beta \gen \sigma \in H_n(M)$. 
    This product is associative, $\SG n$-equivariant, and of degree $n d$, which can be proven exactly as
    the corresponding statements about the product of the rings $H\{\SG n\}$, which are defined in~\cite{lehn-sorger:k3}.
    The product becomes (graded) commutative when restricted to $\Hilb n H$. Thus we have made $\Hilb n H$ a
    graded commutative, unital algebra of degree $n d$.
    \begin{definition}
    	The algebra $\Hilb n H$ is the \emph{$n$-th Hilbert algebra of $H$}.
    \end{definition}
    In case $G$ is trivial, the $n$-Hilbert algebra of $H$ defined here is exactly the algebra
    $\Hilb n H$ of~\cite{lehn-sorger:k3}. For non-trivial $G$, this is no longer true.
    
    The underlying graded vector space of $\sum_{n \ge 0} \Hilb n H(L)$ is naturally isomorphic
    to $S(L) := S^*(\bigoplus_{n \ge 1} H(n L))$, namely as follows:
    Firstly, we introduce linear maps $H_n(L) \to S(L)$, which are defined
    by mapping an element of the form $\sum_{\sigma \in \SG n}
    \bigotimes_{B \in \sigma\backslash[n]} \alpha_{\sigma, B} \gen \sigma$
    to $\frac 1 {n!} \sum_{\sigma \in \SG n} \prod_{B \in \sigma\backslash[n]} \alpha_{\sigma, B}$.
    The restrictions of these morphisms to the $\SG *$-invariant parts define a
    linear map $\bigoplus_{n \ge 0} \Hilb n H(L) \to S(L)$. This map is an isomorphism, which can be
    proven exactly as it is in~\cite{lehn-sorger:k3} for trivial $G$.
    
    Recall that $H^*(X, \set C[2])$ is a (trivially weighted) graded non-counital Frobenius algebra of degree $d$.
    \begin{lemma}
        \label{lem:ring}
        There is a natural isomorphism $\Hilb n {H^*(X, \set C[2])} \to H^*(\Hilb n X, \set C[2n])$ of graded
        unital algebras of degree $nd$.
    \end{lemma}
    
    \begin{proof}
        Recall the just defined isomorphism between the spaces $\bigoplus_{n \ge 0} \Hilb n {H^*(X, \set C[2])}$
        and $S^*(\bigoplus_{\nu > 0} H^*(X, \set C[2]))$ (for $L = \set C$).
        The composition of this isomorphism with
        isomorphism
        between $S^*(\bigoplus_{\nu > 0} H^*(X, \set C[2]))$ and $\bigoplus_{n \ge 0} H^*(\Hilb n X, \set C[2n])$
        of Theorem~\ref{thm:vector_space} induces by restriction the claimed isomorphism of the Lemma on the level
        of graded vector spaces.
        
        That this isomorphism is in fact an isomorphism of unital algebras, is proven in~\cite{lehn-sorger:k3} 
        for $X$ being projective. The proof there does not use the fact that $H^*(X, \set C[2])$ has a counit, in fact
        it only uses its diagonal map.
        It relies on the earlier work
    	in~\cite{lehn:tautological}, which has been extended to the non-projective case above,
        and~\cite{li-qin-wang:vertex}, which can similarly be extended.
        Thus the proof in~\cite{lehn-sorger:k3} also works in the non-projective
    	case, when we replace the notion of a Frobenius algebra by the notion of a non-counital Frobenius algebra.
    \end{proof}
    
    We will now deduce Theorem~\ref{thm:ring} from Lemma~\ref{lem:ring}:
    \begin{proof}[Proof of Theorem~\ref{thm:ring}]
        Let $L, M \in G$.
    
        Let $\lambda = (\lambda_1, \ldots, \lambda_l)$ be a partition of $n$. Let $\nu_i$ the multiplicity of $i$
        in $\lambda$, i.e.~$\lambda = \sum_i \nu_i \cdot i$. Set $\Symm \lambda X := \prod_i \Symm {\nu_i} X$,
        and $\Symm \lambda L := \prod_i \pr_i^* \Symm {\nu_i} L$, where the $\pr_i$ denote the projections onto the
        factors $\Symm{\nu_i} X$.
        Let $\alpha = \sum \alpha_{(1)} \cdots \alpha_{(r)} \in H^*(\Symm \lambda X, \Symm \lambda L[2l])
        = \bigotimes_i S^{\nu_i} H^*(X, L[2])$.
        
        We set
        \[
            \ket \alpha := \sum q(\alpha_{(1)}) \cdots q(\alpha_{(r)})\vac.
        \]
        By Theorem~\ref{thm:vector_space}, the
        cohomology space $H^*(\Hilb n X, \Hilb n L[2n])$ is linearly spanned by classes of the form $\ket \alpha$.
        
        Let $\mu = (\mu_1, \ldots, \mu_m)$ be another partition of $n$
        and $\beta \in H^*(\Hilb \mu X, \Hilb \mu M[2m])$. In order to describe the ring structure of
        $H^*(\Hilb n X, \Hilb n L[2n])$, we have to calculate the classes
        $\ket{\alpha \cup \beta} := \ket \alpha \cup \ket \beta$ in terms of the vector space description given
        by Theorem~\ref{thm:vector_space}.
        
        This means that we have to calculate the numbers
        \[
            \bracket{\gamma \mid \alpha \cup \beta}
            := q(\gamma) \ket{\alpha \cup \beta} \in H^*(\Hilb 0 X, \set C) = \set C
        \]
        for all $\gamma \in H^*_c(\Symm \kappa X, \Symm \kappa {((LM)^{-1})}[2 k])$ for all partitions
        $\kappa = (\kappa_1, \ldots, \kappa_k)$ of $n$,
        and we have to show that they are equal to the numbers that would 
        come out if we calculated the product of $\alpha$ and $\beta$ by the right hand side of the claimed isomorphism
        of the Theorem.

        The class $\ket\alpha$ is given by applying a sequence of correspondences to the vacuum vector: Recall
        from~\cite{nakajima} how to compose correspondences. It follows that $\ket \alpha$ is given by
        \[
            \PD^{-1}(\pr_1)_*(\pr_2^*\alpha \cap \zeta_\lambda),
        \]
        where the symbols have the following meaning:
        The maps $\pr_1$  and $\pr_2$ are
        the projections of $\Hilb n X \times \Symm \lambda X$ onto its factors $\Hilb n X$ and $\Symm \lambda X$. Further,
        $\zeta_\lambda$ is a certain class in $H_*^\BM(Z_\lambda)$,
        where $Z_\lambda$ is the incidence variety
        \[
            Z_\lambda := \Set{(\xi, (\underline x_1, \underline x_2, \ldots)) \in \Hilb n X \times \Symm \lambda X 
            \mid \supp \xi = \sum_i i \underline{x}_i}
        \]
        in $\Hilb n X \times \Symm \lambda X$.
        (Note that $\pr_1^*{\Hilb n L}|_{Z_\lambda} = \pr_2^*{\Symm \lambda L}|_{Z_\lambda}$, and that $p|_{Z_\lambda}$ is
        proper.)
        
        For $\ket\beta$ and $\ket \gamma$ we get similar expressions.
        By definition of the cup-product (pull-back along the diagonal),
        it follows that $\bracket{\gamma \mid \alpha \cup \beta} = \scp{r^* \gamma \cup p^* \alpha \cup q^* \beta}
        {\zeta_{\lambda,
        \mu, \kappa}}$, where $p$, $q$, and $r$ are the projections from $\Symm \lambda X \times \Symm \mu X \times \Symm 
        \kappa X$ onto its three factors,
        and $\zeta_{\lambda, \mu, \kappa}$ is a certain class in $H_*^\BM(Z_{\lambda, \mu, \kappa})$ with
        \[
            Z_{\lambda, \mu, \kappa} := \Set{((\underline x_1, \underline x_2, \ldots), (\underline y_1, 
            \underline y_2, \ldots), (\underline z_1, \underline z_2, \ldots))\mid
                \sum_i i \underline x_i = \sum_j j \underline y_j = \sum_k k \underline z_k}.
        \]
        (The incidence variety is proper over any of the three factors, so everything is well-defined.)        
        The main point is now
        that the incidence variety $Z_{\lambda, \mu, \kappa}$ and the homology class $\zeta_{\lambda, \mu, \kappa}$
        are independent of the locally constant systems $L$ and $M$.
        In particular, we can calculate $\zeta_{\lambda, \mu, \kappa}$
        once we know the cup-product in the case $L = M = \set C$. But this is the case that is described
        in Lemma~\ref{lem:ring}, which we will analyse now.
        
        First of all, the incidence variety is given by
        \[
            Z_{\lambda, \mu, \kappa} = \sum_{\sigma, \tau} Z_{\sigma, \tau}  
        \]
        where $\sigma$ and $\tau$ run through all permutations with cycle type $\lambda$ and $\mu$, respectively,
        such that
        $\rho := \sigma\tau$ has cycle type $\kappa$. The
        varieties $Z_{\sigma, \tau}$ are defined as follows:
        
        As the orbits of the group action of $\langle\sigma\rangle$ on $[n]$ correspond to the entries of
        the partition $\lambda$,
        there exists a natural map $X^{\sigma\backslash[n]} \to \Symm \lambda X$, which is given by symmetrising.
        Furthermore
        the natural surjection $\sigma\backslash[n] \to \langle\sigma, \tau\rangle\backslash[n]$ induces
        a diagonal embedding
        $X^{\langle\sigma, \tau\rangle\backslash[n]} \to X^{\sigma\backslash[n]}$. Composing
        both maps, we get a
        natural map $X^{\langle\sigma, \tau\rangle\backslash[n]} \to \Symm \lambda X$. Analoguously, we get maps from
        $X^{\langle\sigma, \tau\rangle\backslash[n]}$
        to $\Symm \mu X$ and $\Symm \kappa X$. Together, these maps define a diagonal embedding
        \[
            i_{\tau, \sigma}\colon
            X^{\langle\sigma, \tau\rangle\backslash[n]} \to \Symm \kappa X \times \Symm \lambda X \times \Symm \mu X.
        \]
        We define $Z_{\sigma, \tau}$ to be the image of this map.
        
        By Lemma~\ref{lem:ring}, the class $\zeta_{\lambda, \mu, \kappa}$ is given by
        $\sum_{\sigma, \tau} (i_{\sigma, \tau})_* \zeta_{\sigma, \tau}$,
        where each class $\zeta_{\sigma, \tau} \in H^\BM_*(X^{\langle\sigma, \tau\rangle\backslash [n]})$
        is Poincar\'e dual to
        $c_{\sigma, \tau} \, e^{\gamma(\sigma, \tau)}$. Here, $c_{\sigma, \tau}$ is a certain
        combinatorial factor (possibly depending on $\sigma$ and $\tau$), whose precise value is of no
        concern for us.
        
        Having derived the value of $\zeta_{\lambda, \mu, \kappa}$ from Lemma~\ref{lem:ring}, we have thus
        calculated the value $\bracket{\gamma \mid \alpha \cup \beta}$.
        
        Now we have to compare this value with the one that is predicted
        by the description of the cup-product given by the right hand side of the claimed isomorphism of the Theorem.
        With the same analysis as above, we find this value is also given by a correspondence on $Z_{\lambda, \mu, \kappa}$
        with the class $\sum_{\tau, \sigma} (i_{\sigma, \tau})_* c_{\sigma, \tau} \PD(e^{\gamma(\sigma, \tau)})$
        with the same combinatorial factors $c_{\sigma, \tau}$ as above. We thus find that the claimed ring structure
        yields the correct value of $\bracket{\gamma \mid \alpha \cup \beta}$.
    \end{proof}
    
    \begin{remark}
        One can also define a natural diagonal map for the Hilbert algebras $\Hilb n H$ making them into graded,
        non-counital Frobenius algebras of degree $nd$. The isomorphism of Theorem~\ref{thm:ring} then becomes an
        isomorphism of graded non-counital Frobenius algebras.
    \end{remark}
    
    \section{The generalised Kummer varieties}
    
    Finally, we want to use Theorem~\ref{thm:ring} to study the cohomology ring of the generalised Kummer varieties.
    
    Let $H$ be a non-counital graded Frobenius algebra of degree $d$ that is moreover endowed with a compatible
    structure of a
    cocommutative Hopf algebra of degree $d$. The
    comultiplication $\delta$ of the Hopf algebra structure is of degree $-d$. The counit of the Hopf algebra structure 
    is denoted by $\epsilon$ and is of degree $d$.
    We further assume that $H$ is also equipped with a $G$-weighting for a finite group $G$.
    \begin{example}
        Let $X$ be an abelian surface. The group structure on $X$ induces naturally a graded Hopf algebra structure
        of degree $2$
        on the graded Frobenius algebra $H^*(X, \set C[2])$. This algebra is also trivially $X[n]\dual$-weighted,
        where $X[n]\dual \simeq (\set Z/(n))^4$
        is the character group of the group of $n$-torsion points on $X$. (Trivially weighted means
        that the only non-trivial $X[n]\dual$-weight space of $H^*(X, \set C[2])$ is the one corresponding to the 
        identity element $0$.)
    \end{example}

    Let $n$ be a positive integer. Recall the definition of the ($G$-weighted) Hilbert algebra $\Hilb n H$.
    Repeated application of the comultiplication $\delta$
    induces a map $\delta\colon H \to H^{\otimes n} = H^{\id\backslash [n]}$, which is of degree $-(n - 1) d$.
    Its image lies in the subspace of symmetric tensors. Thus we can define
    a map $\phi\colon H \to \Hilb n H$
    with $\phi(\alpha) := \delta(\alpha) \gen \id$.
    One can easily check that this map is an algebra homomorphism of degree $-(n - 1) d$, making $\Hilb n H$
    into an
    $H$-algebra.
    
    Define
    \[
        \Kummer n H := \Hilb n H \otimes_H \set C,
    \]
    where we view $\set C$ as an $H$-algebra of
    degree $d$ via the Hopf counit $\epsilon$. It is $\Kummer n H$ a ($G$-weighted)
    graded Frobenius algebra of degree $nd$.
    \begin{definition}
        The algebra $\Kummer n H$ is the \emph{$n$-th Kummer algebra of $H$}.
    \end{definition}
    The reason of this naming is of course Theorem~\ref{thm:kummer}.
    
    \begin{proof}[Proof of Theorem~\ref{thm:kummer}]
        Let $n\colon X \to X$ denote the morphism that maps $x$ to $n \cdot x$.
        There is a natural cartesian square
        \begin{equation}
            \label{eq:kummer}
            \begin{CD}
                X \times \Kummer n X @>{\nu}>> \Hilb n X \\
                @VpVV @VV{\sigma}V \\
                X @>>n> X,
            \end{CD}
        \end{equation}
        where $p$ is the projection on the first factor and $\nu$ maps a pair $(x, \xi)$ to $x + \xi$, the subscheme
        that is given by translating $\xi$ by $x$ (\cite{beauville:kummer}). 
        Let $G$ be the character group of the Galois group of $n$, i.e.~$G = X[n]\dual$.
        Each element $L$ of $G$ corresponds to a locally
        constant system $L$ on $X$, and we have $n_* \set C = \bigoplus_{L \in G} L$.
        It follows that $\nu$ is the $G$-covering of $\Hilb n X$ with $\nu_* \set C = \bigoplus_{L \in G} \Hilb n L$.
        
        Together with Theorem~\ref{thm:ring}, this leads to the claimed description of the cohomology ring
        of $\Kummer n X$: Firstly, there is a natural isomorphism
        \[
            H^*(\Kummer n X, \set C[2n]) \to H^*(X \times \Kummer n X, \set C[2n]) \otimes_{H^*(X, \set C[2])} \set C
        \]
        of unital algebras (the tensor product is taken with
        respect to the map $p^*$ and the Hopf counit $H^*(X, \set C[2]) \to \set C$). By the Leray spectral sequence
        for $\nu$ and by~\eqref{eq:kummer}, the right hand side is naturally isomorphic to
        \[
            H^*(\Hilb n X, \nu_* \set C[2n]) \otimes_{H^*(X, \set C[2])} \set C
            = \bigoplus_{L \in G} H^*(\Hilb n X, \Hilb n L[2n]) \otimes_{H^*(X, \set C[2])} \set C
        \]
        (where the tensor
        product is taken with respect to the map $\sigma^*$ and the Hopf counit).
        
        By Theorem~\ref{thm:ring}, the algebra $\bigoplus_{L \in G} H^*(\Hilb n X, \Hilb n L[2n])$ is naturally 
        isomorphic to $\Hilb n {\bigoplus_{L \in G} H^*(X, L[2])}$. Now $H^*(X, L[2]) = 0$ unless $L$ is the
        trivial bundle, which follows from the fact that all classes in $H^*(X, \set C)$ are invariant under the action of 
        the Galois group of $n$, i.e.~correspond to the trivial character.
        Thus there is a natural isomorphism
        \[
            \bigoplus_{L \in G} H^*(\Hilb n X, \Hilb n L[2n]) \to \Hilb n {H^*(X, \set C[2])},
        \]
        of $G$-weighted algebras, where we endow $H^*(X, \set C[2])$ with the trivial $G$-weighting.
        Under this isomorphism, the map $\sigma^*$ corresponds
        to the homomorphism $\phi$ defined above by Example~\ref{ex:hilbert-chow}. Thus
        we have proven the existence of a natural isomorphism
        \[
            H^*(\Kummer n X, \set C[2n]) \to \Hilb n {H^*(X, \set C[2])} \otimes_{H^*(X, \set C[2])} \set C
        \]
        of unital, graded algebras.
        But the right hand side is nothing but $\Kummer n H$, thus the Theorem is proven.
    \end{proof}

	\appendix
	\bibliographystyle{amsalpha}
	\bibliography{my}
	
\end{document}